\documentclass[12pt]{article}
\usepackage{mathtools}
\usepackage[utf8]{inputenc}
\usepackage{amssymb,amsmath,amsfonts,eucal,mathrsfs,amsthm} 
\usepackage[colorinlistoftodos]{todonotes}
\usepackage{xcolor}
\usepackage{enumitem}
\usepackage[notcite,notref]{showkeys}
\newtheorem{theorem}{Theorem}
\newtheorem{proposition}[theorem]{Proposition}
\newtheorem{lemma}[theorem]{Lemma}

\theoremstyle{definition}
\newcommand{\R}{\mathbb{R}}
\newcommand{\Z}{\mathbb{Z}}
\newcommand{\Q}{\mathbb{Q}}
\newcommand{\Sf}{\mathbb{S}}

\newcommand{\spa}{\mbox{span}}

\newcommand{\Ric}{\mbox{Ric}}

\newcommand{\trace}{\mbox{tr\,}}

\def\<{{\langle}}
\def\>{{\rangle}}
\def\B{\mathcal{B}}
\def\CP{\mathord{\mathbb C}\mathord{\mathbb P}}

\def\n{\nabla}

\def\a{\alpha}

\def\be{\begin{equation} }
\def\ee{\end{equation} }

\newcommand\blfootnote[1]{\begingroup
\renewcommand\thefootnote{}\footnote{#1}
\addtocounter{footnote}{-1}
\endgroup}
\begin{document}

\title{Ricci pinched compact submanifolds\\ in space forms}
\author{M. Dajczer and Th. Vlachos}
\date{}
\maketitle

\begin{abstract}
We investigate the compact submanifolds in Riemannian space 
forms of nonnegative sectional curvature that satisfy a lower 
bound on the Ricci curvature, that bound depending solely on the 
length of the mean curvature vector of the immersion. 
While generalizing the results, we give a 
positive answer to a conjecture by H. Xu and J. Gu in 
(2013, Geom. Funct. Anal. 23). Our main accomplishment 
is the elimination of the need for the mean curvature 
vector field to be parallel.
\end{abstract}
\blfootnote{\textup{2020} \textit{Mathematics Subject 
Classification}: 53C20, 53C40, 53C42.}

Let $f\colon M^n\to\Q^{n+m}_c$ be an isometric immersion 
of a compact Riemannian manifold of dimension $n$ into
a simply-connected space form of sectional curvature $c$   
and substantial codimension $m$. Throughout this paper, 
the (not normalized) Ricci curvature of $M^n$ is assumed 
to satisfy at any point the pinching condition
\be
\Ric_M\geq (n-2)(c+H^2)\tag{$\ast$},
\ee 
where $H$ stands for the norm of the (normalized) mean 
curvature vector. 

By the condition $(*)$ being satisfied with equality at 
$x\in M^n$ we mean that the inequality at that point
is not strict, that is, there exists a unit vector 
$X\in T_xM$ such that $\Ric_M(X)=(n-2)(c+H^2)$. If 
otherwise, we call the inequality $(*)$ at $x\in M^n$ strict.
Notice that if the inequality holds strictly at any given point, 
it will persist in its strict form after subjecting the 
submanifold to a sufficiently small smooth deformation.
\vspace{1ex}

A complete classification of the submanifolds as above for $c>0$
was obtained by N. Ejiri \cite{E} in 1979 under the assumptions 
that $f$ is a minimal immersion and that the manifold $M^n$ is 
both, oriented and simply connected. 
H. Xu and J. Gu \cite{XG} in 2013 generalized  Ejiri's result 
under the assumption that $c+H^2>0$ by only requiring 
the mean curvature vector to be parallel and the submanifold  
orientable. Additionally, they proved that when the assumption 
regarding the mean curvature is removed, but a strict inequality 
in $(*)$ holds at every point, the manifold must be homeomorphic 
to a sphere.

In the same paper, Xu and Gu put forward a conjecture, suggesting 
that their findings should hold true even when eliminating the 
condition on the mean curvature vector altogether. 
In this paper, we give affirmative confirmation to their conjecture. 
To achieve that result we do not require the manifold to be oriented. 
Furthermore, to reach the conclusion that the manifold is 
homeomorphic to a sphere we only ask for strict inequality 
in the condition $(*)$ at some point.

\begin{theorem}\label{thm1}
Let $f\colon M^n\to\mathbb{Q}_c^{n+m}$, $n\geq 4, c\geq 0$, 
be an isometric immersion of a compact Riemannian manifold  
that satisfies  the condition $(*)$ at any point. Then 
either $M^n$ is homeomorphic to $\Sf^n$ or one of the 
following holds:
\begin{itemize}
\item[(i)] $M^n$ is 
$\Sf^{n/2}(r/\sqrt{2})\times\Sf^{n/2}(r/\sqrt{2})$, 
$r=1/\sqrt{c+H^2}$, and $f=j\circ g$, 
where $g\colon M^n\to\Sf^{n+1}(r)$ is the standard embedding 
and $j\colon\Sf^{n+1}(r)\to\Q_c^{n+2}$ an umbilical inclusion.
\item[(ii)] $M^n$ is the projective
plane $\CP_r^2$ of constant holomorphic curvature $4r^2/3$ with
$r=1/\sqrt{c+H^2}$ and $f=j\circ g$, where
$g$ is the standard immersion of  $\CP_r^2$ into $\Sf^7(r)$
and $j\colon\Sf^7(r)\to\Q_c^8$ an umbilical inclusion.
\end{itemize}
\end{theorem}

\noindent\emph{Remarks}: $(i)$ In both cases 
the umbilical inclusion may be totally geodesic.

\noindent $(ii)$ If $n$ is even and $c>0$ the result 
generalizes Theorem $2$ in \cite{DV} by means of a lower 
bound for the Ricci curvature.
\vspace{1ex}

When $M^n$ possesses the topological structure of a 
sphere, the conjecture put forth by Xu and Gu proposes 
that it should not merely be topologically equivalent but 
diffeomorphic to a sphere. This holds true for dimensions 
$n=5,6,12,56,61$ as in these cases, it is established that 
the differentiable structure is unique; as stated by 
Corollary $1.15$ in \cite{WX}.

\section{The pinching condition}

In this section, we analyze the relation between our pinching 
condition $(*)$  and the one for $c>0$ due to Lawson and Simons   
in their seminal paper \cite{LS}. Their result, was later 
strengthened by Elworthy and Rosenberg \cite[p.\ 71]{ER} by  
only requiring the bound to be strict at some point of the 
submanifold. The case $c=0$ was later considered by Xin \cite{Xi}.

\begin{theorem}{\em(\cite{ER},\cite{LS},\cite{Xi})}\label{ls} 
Let $f\colon M^n\to\Q_c^{n+m}$, $n\geq 4,c\geq 0$, 
be an isometric immersion of a compact manifold and  $p$ an integer
such that $1\leq p\leq n-1$. Assume that at any point $x\in M^n$ 
and for \emph{any} orthonormal basis $\{e_1,\ldots,e_n\}$ of $T_xM$ 
the second fundamental form $\alpha_f\colon TM\times TM\to N_fM$  
satisfies
\be
\Theta_p=\sum_{i=1}^p\sum_{j=p+1}^n\big(2\|\a_f(e_i,e_j)\|^2
-\<\a_f(e_i,e_i),\a_f(e_j,e_j)\>\big)\leq p(n-p)c \tag{$\#$}.
\ee
If the inequality $(\#)$ is strict at a point of $M^n$, then 
there are no stable $p$-currents and the homology groups satisfy 
$H_p(M^n;\mathbb{Z})=H_{n-p}(M^n;\mathbb{Z})=0$.
\end{theorem}

Recall that a vector in the normal bundle 
$\eta\in N_fM(x)$  at $x\in M^n$ is named a 
\emph{Dupin principal normal} of $f\colon M^n\to\Q^{n+m}_c$ 
if the vector subspace
$$
E_\eta(x)=\left\{X\in T_xM\colon\alpha_f(X,Y)
=\<X,Y\>\eta\;\,\text{for all}\;\,Y\in T_xM\right\}
$$
is at least two dimensional. That dimension is called the 
\emph{multiplicity} of $\eta$.
\vspace{1ex}

The proof of the following results is inspired by computations
given by us in \cite{DV} and by Xu and Gu in \cite{XG}.

\begin{proposition}\label{prop}
Let $f\colon M^n\to\Q_c^{n+m},n\geq4$, be an isometric 
immersion satisfying the inequality $(*)$ at $x\in M^n$.
Then the following assertions at $x\in M^n$ hold:
\vspace{1ex}

\noindent $(i)$ The inequality $(\#)$ is satisfied for any 
integer $2\leq p\leq n/2$ and any orthonormal basis of $T_xM$. 
Moreover, if the inequality $(*)$ is strict then 
also $(\#)$ is strict for any integer $2\leq p\leq n/2$. 
\vspace{1ex}

\noindent $(ii)$ Assume that equality holds in $(\#)$ for 
a certain integer $2\leq p\leq n/2$ and a given orthonormal 
basis $\{e_j\}_{1\leq j\leq n}$ of $T_xM$. Then the Ricci 
tensor satisfies 
$$
\Ric_M(X)=(n-2)(c+H^2)\;\; \text{for any unit}\;\; X\in T_xM.
$$
Moreover, we have:
\begin{itemize}
\item[(a)] If $n\geq 5$ then  either $c=0$ and $f$ is totally 
geodesic or we have that $p=n/2$ and there are distinct Dupin
principal normals $\eta_1$ and $\eta_2$ such that 
$E_{\eta_1}=\spa\{e_1,\dots,e_p\}$ and 
$E_{\eta_2}=\spa\{e_{p+1},\dots,e_n\}$.
\item[(b)] If $n=4$ and $p=2$ there are normal vectors $\eta_j$, 
$j=1,2$, such that 
\be\label{Vj}
\pi_{V_j}\circ A_\xi|_{V_j}=\<\xi,\eta_j\>I\;\,
\text{for any}\;\,\xi\in N_f(x)
\ee
where $V_1=\spa\{e_1,e_2\}$, $V_2=\spa\{e_3,e_e\}$ and 
$\pi_{V_j}\colon T_xM\to V_j$ denotes the projection.
\end{itemize}
\end{proposition}

\proof We first recall that the Gauss equation of 
$f\colon M^n\to\Q_c^{n+m}$ implies that the Ricci curvature  
for any unit vector $X\in T_xM$ satisfies 
\be\label{ric}
\Ric_M(X)=(n-1)c+\sum_{\a=1}^m(\mathrm{tr}A_\a)
\<A_\a X,X\>-\sum_{\a=1}^m\|A_{\a}X\|^2, 
\ee
where the $A_{\a}$, $1\leq\a\leq m$, are the shape operators 
of $f$ associated to an orthonormal basis 
$\{\xi_\alpha\}_{1\leq\alpha\leq m}$ of the normal space 
$N_fM(x)$ of the submanifold.

From now on, we agree that  
$\{\xi_\alpha\}_{1\leq\alpha\leq m}$ satisfies that the 
(normalized) mean curvature vector is $\mathcal{H}(x)=H(x)\xi_1$ 
when $H(x)\neq 0$. For a given orthonormal basis 
$\{e_j\}_{1\leq j\leq n}$ of $T_xM$, we denote for simplicity 
$\a_{ij}=\alpha_f(e_i,e_j)$, $1\leq i,j\leq n$. Then, we have
\begin{align}\label{a}
\Theta_p&=\;2\sum_{i=1}^{p}\sum_{j=p+1}^n\|\a_{ij}\|^2-n
\sum_{i=1}^p\<\a_{ii},\mathcal{H}\>
+\sum_{i,j=1}^p\<\a_{ii},\a_{jj}\>\nonumber\\
=&\;2\sum_{i=1}^p\sum_{j=p+1}^n\sum_\a\<A_{\alpha }e_i,e_j\>^2-nH
\sum_{i=1}^p\<A_1e_i,e_i\>
+\sum_\a\big(\sum_{i=1}^p\<A_\alpha e_i,e_i\>\big)^2\nonumber\\
\leq&\; 2\sum_{i=1}^p\sum_{j=p+1}^n\sum_\a\<A_{\alpha }e_i,e_j\>^2-nH
\sum_{i=1}^p\<A_1e_i,e_i\>
+p\sum_\a\sum_{i=1}^p\<A_{\alpha}e_i,e_i\>^2,
\end{align}
where the inequality part was obtained using the Cauchy-Schwarz 
inequality
\be\label{a-}
\big(\sum_{i=1}^p\<A_\alpha e_i,e_i\> 
\big)^2\leq p\sum_{i=1}^p\<A_{\alpha}e_i,e_i\> ^2. 
\ee
Since $p\geq 2$ by assumption, then 
\begin{align}\label{b-}
2\sum_{i=1}^p&\sum_{j=p+1}^n\sum_\a\<A_{\alpha}e_i,e_j\>^2
+p\sum_{i=1}^p\sum_\a\<A_{\alpha}e_i,e_i\>^2\nonumber\\
&\leq p\sum_{i=1}^p\sum_{j=p+1}^n
\sum_\a\<A_{\alpha }e_i,e_j\>^{2}+p\sum_{i=1}^p\sum_\a
\<A_{\alpha}e_i,e_i\>^2\nonumber\\
&\leq p\sum_{i=1}^p\sum_\a\|A_{\alpha}e_i\|^2 
\end{align}
and thus \eqref{a} implies that
$$
\Theta_p\leq p\sum_{i=1}^p\sum_\a\|A_{\alpha}e_i\|^2
-nH\sum_{i=1}^p\<A_1e_i,e_i\>.
$$
Setting $\varphi=A_1-HI$ and using \eqref{ric}, we obtain 
\begin{align}\label{first}
\Theta_p
&\leq p\sum_{i=1}^p((n-1)c-\Ric_M(e_i))
+(p-1)nH\sum_{i=1}^p\<A_1e_i,e_i\>\nonumber\\
&= p\sum_{i=1}^p((n-1)(c+H^2)-\Ric_M(e_i))
-p(n-p)H^2\nonumber\\
&\;+(p-1)nH\sum_{i=1}^p\<\varphi e_i,e_i\>.
\end{align}
Then, we have
\begin{align}\label{d}
\Theta_p
&\leq p^2\big((n-1)(c+H^2)-\Ric_M^{\text{min}}(x)\big)-p(n-p)H^2\nonumber\\
&\;+(p-1)nH\sum_{i=1}^p\<\varphi e_i,e_i\>.
\end{align}
where 
$\Ric_M^{\text{min}}(x)=\min\left\{\Ric_M(X)\colon X\in T_xM,\|X\|=1\right\}$.
\vspace{1ex}

We claim that 
\be\label{min}
(n-1)(c+H^2)\geq\Ric_M^{\text{min}}(x)
\ee
and that equality holds if $f$ is umbilical at $x\in M^n$.
Indeed, it follows from the Gauss equation that the scalar 
curvature $\tau$ is given by
\be\label{tau}
\tau=n(n-1)c-S+n^2H^2,
\ee
where $S$ denotes the norm of the second fundamental form. 
Hence 
\be\label{sca}
S\leq n(n-1)c+n^2H^2-n\Ric_M^{\text{min}}(x).
\ee
Therefore, we have
\begin{align*}
(n-1)(c+H^2)-&\Ric_M^{\text{min}}(x)\geq\frac{1}{n}(S-nH^2)
=\frac{1}{n}\|\phi\|^2,
\end{align*}
where $\phi=\a-\<\,,\,\>\mathcal H$ is the traceless second 
fundamental form, and the claim now follows.

From \eqref{min} and having that $p\leq n/2$, then
\be\label{re}
p^2\big((n-1)(c+H^2)-\Ric_M^{\text{min}}(x)\big)\leq p(n-p)
\big((n-1)(c+H^2)-\Ric_M^{\text{min}}(x)\big).
\ee
Therefore, it follows from \eqref{d} the estimate
\be\label{c1}
\Theta_p
\leq p(n-p)\big((n-1)(c+H^2)-\Ric_M^{\text{min}}(x)-H^2\big)
+(p-1)nH\sum_{i=1}^p\<\varphi e_i,e_i\>.
\ee

Next, we provide a second estimate of 
$$
\Theta_p=\;\sum_\a\Big(2\sum_{i=1}^p\sum_{j=p+1}^n\<A_\a e_i,e_j\>^2
-\sum_{i=1}^p\<A_\a e_i,e_i\>\sum_{j=p+1}^n\<A_\a e_j,e_j\>\Big)
$$
Decomposing
\begin{align*}
&\sum_{i=1}^p\<A_\a e_i,e_i\>\sum_{j=p+1}^n\<A_\a e_j,e_j\>\\
=&\frac{n-p}{n}\sum_{i=1}^p\<A_\a e_i,e_i\>
\sum_{j=p+1}^n\<A_\a e_j,e_j\>
+  \frac{p}{n}\sum_{i=1}^p\<A_\a e_i,e_i\>
\sum_{j=p+1}^n\<A_\a e_j,e_j\>,
\end{align*}
we have
\begin{align*}
&\Theta_p
=\sum_\a\Big(2\sum_{i=1}^p\sum_{j=p+1}^n\<A_\a e_i,e_j\>^2
-\frac{n-p}{n}\trace A_\a\sum_{i=1}^p\<A_\a e_i,e_i\>\\
&+\frac{n-p}{n}\big(\sum_{i=1}^p\<A_\a e_i,e_i\>\big)^2
-\frac{p}{n}\trace A_\a\sum_{j=p+1}^n\<A_\a e_j,e_j\>
+\frac{p}{n}\big(\sum_{j=p+1}^n\<A_\a e_j,e_j\>\big)^2\Big).
\end{align*}
Using the Cauchy-Schwarz inequality, we obtain
\begin{align*}\label{a1}
\Theta_p
\leq&\sum_\a\Big(2\sum_{i=1}^p\sum_{j=p+1}^n\<A_\a e_i,e_j\>^2
-\frac{n-p}{n}\trace A_\a\sum_{i=1}^p\<A_\a e_i,e_i\>\\
&+\;\frac{p(n-p)}{n}\sum_{i=1}^p\<A_\a e_i,e_i\>^2
-\frac{p}{n}\trace A_\a\sum_{j=p+1}^n\<A_\a e_j,e_j\>\\
&+\;\frac{p(n-p)}{n}\sum_{j=p+1}^n\<A_\a e_j,e_j\>^2\Big)\\
\leq&\;\frac{p(n-p)}{n}S-(n-p)H\sum_{i=1}^p\<A_1e_i,e_i\>
-pH\sum_{j=p+1}^n\<A_1e_j,e_j\>\nonumber\\
=&\;\frac{p(n-p)}{n}S-pnH^2-(n-2p)H\sum_{i=1}^p\<A_1e_i,e_i\>\\
=&\;\frac{p(n-p)}{n}S-2p(n-p)H^2-(n-2p)H
\sum_{i=1}^p\<\varphi e_i,e_i\>.
\end{align*}
Then using \eqref{sca} we obtain
\be\label{a2}
\Theta_p
\leq p(n-p)\big((n-1)(c+H^2)-\Ric_M^{\text{min}}(x)-H^2\big)
-(n-2p)H\sum_{i=1}^p\<\varphi e_i,e_i\>.
\ee

Finally, by computing $(n-2p)\times$\eqref{c1}$+n(p-1)\times$\eqref{a2} 
and using $(*)$ it follows that
\be\label{ac1}
\Theta_p
\leq p(n-p)\Big((n-1)(c+H^2)-\Ric_M^{\text{min}}(x)-H^2\Big)\nonumber\\
\leq p(n-p)c,
\ee
and $(\#)$ has been proved. Clearly, if the inequality $(*)$ is 
strict then also \eqref{ac1} becomes strict, and this completes 
the proof of part $(i)$.
\vspace{1ex}

We now prove part $(ii)$. Thus, we assume that 
equality holds in $(\#)$  for a certain integer $2\leq p\leq n/2$ 
and a given orthonormal basis $\{e_j\}_{1\leq j\leq n}$ of $T_xM$. 
Then, all the inequalities from \eqref{a} to \eqref{d} as well as 
the ones from \eqref{re} to \eqref{ac1} become equalities. 
In particular, from \eqref{a-} we obtain that
\be\label{r}
\<A_\a e_i,e_i\>=\rho_\a\;\,\text{for all}
\;1\leq i\leq p,\;\,1\leq\a\leq m.
\ee
From \eqref{b-} it follows that
\be\label{k}
(p-2)\<A_\a e_i,e_j\>=0\;\,\text{for all}
\;\,1\leq i\leq p,\;p+1\leq j\leq n,\;1\leq\a\leq m, 
\ee
and
\be\label{ii}
\<A_\a e_i,e_{i'}\>=0\;\,\text{for all}\;\,1\leq i\neq i' 
\leq p,\;1\leq\a\leq m.
\ee
We obtain from \eqref{re} that 
\be\label{p-min}
(p-n/2)\big(\Ric_M^{\text{min}}(x)-(n-1)(c+H^2) \big)=0.  
\ee
From \eqref{first} and \eqref{d} we have 
$\Ric_M(e_i)=\Ric_M^{\text{min}}(x)$ and then \eqref{ac1} gives 
\be\label{ricb}
\Ric_M(e_i)=\Ric_M^{\text{min}}(x)=(n-2)(c+H^2)
\;\,\text{for all}\;\,1\leq i\leq p.
\ee

At first suppose that $p\neq n/2$. Then \eqref{p-min} implies 
that equality holds in the inequality \eqref{min}, and we have
seen that this gives that $f$ is umbilical at $x$.
At umbilical points, it now follows from \eqref{ricb} that 
$\Ric_M(x)=(n-1)(c+H^2)$, that $c=0$ and $f$ is 
totally geodesic. 

Hereafter, let $2p=n$. Then equality also holds in $(\#)$ for
the reordered orthonormal basis 
$\{e_{p+1},\dots,e_n, e_1,\dots,e_p\}$ of $T_xM$. 
Thus, we also have 
\be\label{r1}
\<A_\a e_j,e_j\>=\mu_\a\;\,\text{for all}
\;\,p+1\leq j\leq n,\;\,1\leq\a\leq m,
\ee
\be\label{jj}
\<A_\a e_j,e_{j'}\>=0\;\,\text{for all}\;\,p+1\leq j\neq j' 
\leq n,\;1\leq\a\leq m,
\ee
and 
\be\label{ricbj}
\Ric_M(e_j)=\Ric_M^{\text{min}}(x)=(n-2)(c+H^2)
\;\,\text{for all}\;\,p+1\leq j\leq n.
\ee
Hence, we obtain from \eqref{ricb} and \eqref{ricbj}  
that the Ricci tensor satisfies 
\be\label{ein}
\Ric_M(X)=(n-2)(c+H^2)\;\,\text{for any unit}\;\,X\in T_xM.
\ee

In particular, if $n\geq 6$ then it follows from \eqref{r}, 
\eqref{k}, \eqref{ii}, \eqref{r1} and \eqref{jj} that 
$\eta_1=\sum_\a\rho_\a\xi_\a$ and 
$\eta_2=\sum_\a\mu_\a\xi_\a$ are Dupin principal normals with 
$$
E_{\eta_1}=\spa\{e_1,\dots,e_p\}\;\;
{\text {and}}\;\;E_{\eta_2}=\spa\{e_{p+1},\dots,e_n\}.
$$
If $\eta_1=\eta_2$, then $f$ is umbilical at $x$ and 
equality holds in \eqref{min}. This combined with \eqref{ein} 
yields $c=0$ and that $f$ is totally geodesic at $x$. 
Otherwise $\eta_1$ and $\eta_2$ are distinct Dupin principal 
normals, and this concludes the proof of part~$(a)$. 

Finally, if $n=4$ then for any $\xi\in N_f(x)$ we have 
\eqref{Vj} where $\eta_1=\sum_\a\rho_\a\xi_\a$ and 
$\eta_2=\sum_\a\mu_\a\xi_\a$, and part $(ii)$ has also 
been proved.\qed

\section{The proof of Theorem \ref{thm1}}

We first establish a topological result 
necessary for proving the theorem.

\begin{lemma}\label{p1}
Let $f\colon M^n\to\Q_c^{n+m}$, $n\geq 4,c\geq0$, 
be an isometric immersion of a compact manifold satisfying
\be\label{p=1}
\Ric_M\geq \frac{n(n-1)}{n+2}(c+H^2) 
\ee
with strict inequality at some point.
Then $\pi_1(M^n)=0$ and $H_{n-1}(M^n,\Z)=0$.
\end{lemma}

\proof From \eqref{tau} and \eqref{p=1} we obtain that
$\phi=\a-\<\,,\,\>\mathcal H$ satisfies 
\be\label{phi}
\|\phi\|^2\leq\frac{2n(n-1)}{n+2}(c+H^2).
\ee
Let $\{e_i\}_{1\leq i\leq n}$ be an orthonormal tangent
basis and let  $\{\xi_\a\}_{1\leq\a\leq m}$ be an orthonormal normal 
basis at $x\in M^n$.
Using \eqref{ric} we have 
\begin{align*}
\sum_{j=2}^n&\big(2\|\a_{1j}\|^2
-\<\a_{11},\a_{jj}\>\big)\\
=&\;2\sum_\a\sum_{j=2}^n\<A_\a e_1,e_j\>^2
-\sum_\a\<A_\a e_1,e_1\>\sum_{j=2}^n\<A_\a e_j,e_j\>\\
=&\;\sum_\a\sum_{j=2}^n\<A_\a e_1,e_j\>^2
-\sum_\a\trace A_\a\<A_\a e_1,e_1\>
+\sum_\a\|A_\a e_1\|^2\\
=&\;\sum_{j=2}^n\|\phi( e_1,e_j)\|^2+(n-1)c-\Ric_M(e_1).
\end{align*}
This inequality together with \eqref{p=1} and \eqref{phi} give
\begin{align*}
\sum_{j=2}^n&\big(2\|\a_{1j}\|^2
-\<\a_{11},\a_{jj}\>\big)
\leq\;\frac{1}{2}\|\phi\|^2+(n-1)c-\Ric_M(e_1)\\
\leq&(n-1)c+ \frac{n(n-1)}{n+2}(c+H^2)-\Ric_M(e_1)
\leq(n-1)c.
\end{align*}
Clearly, if \eqref{p=1} is a strict inequality at a point 
then also is the above. Then by Theorem \ref{ls} there 
are no stable $1$-currents on $M^n$ and thus 
$H_1(M^n,\Z)=H_{n-1}(M^n,\Z)=0$. Since in each nontrivial free 
homotopy class there is a length minimizing curve, we conclude 
that $\pi_1(M^n)=0$.\qed
\vspace{2ex}

\noindent\emph{Proof of Theorem \ref{thm1}:} 
According to part $(i)$ of Proposition \ref{prop} the 
inequality $(\#)$ is satisfied at any point of $M^n$ 
for any $2\leq p\leq n/2$ and any orthonormal tangent 
basis at that point. Since $n(n-1)/(n+2)<(n-2)$, then
for $c>0$ we have that
\be\label{p=1+}
\Ric_M\geq(n-2)(c+H^2)>\frac{n(n-1)}{n+2}(c+H^2).
\ee
If $c=0$, then compactness implies that there exists 
a point $x\in M^n$ where $H(x)\neq 0$, and  
\eqref{p=1+} holds at that point. Now Lemma \ref{p1} 
yields that $M^n$ is simply connected and 
that $H_{n-1}(M^n,\Z)=0$.
\vspace{1ex}

We need to distinguish two cases:
\vspace{1ex}

\noindent\emph{Case I}. Suppose first that 
\be\label{Hom}
H_p(M^n;\Z)=0=H_{n-p}(M^n;\Z)\;\,{\text {for all}}\;\,2\leq p\leq n/2, 
\ee
which by Theorem \ref{ls} is necessarily the case if $(*)$ 
is strict at some point of $M^n$. Hence $M^n$ is a simply 
connected homology sphere and it follows from  the 
Hurewicz isomorphism theorem that $M^n$ is a homotopy 
sphere.  Then the resolution of the generalized Poincar\'{e} 
conjecture gives that $M^n$ is homeomorphic to $\Sf^n$. 
\vspace{1ex}

\noindent\emph{Case II}. Suppose now that \eqref{Hom} does 
not hold. Consider the nonempty set
$$
P=\left\{2\leq p\leq n/2: H_p(M^n;\Z)
\neq0\;\;{\text {or}}\;\;H_{n-p}(M^n;\Z)\neq0\right\}
$$
and set $k=\max P$. Hence $H_k(M^n;\Z)\neq 0$ or 
$H_{n-k}(M^n;\Z)\neq 0$. Then, by Theorem \ref{ls}, there 
exists at any point $x\in M^n$ an orthonormal tangent basis
such that equality holds in $(\#)$ for $p=k$. Moreover, 
at any point $x\in M^n$ we have from part $(ii)$ of 
Proposition \ref{prop} that 
$$
\Ric_M(X)=(n-2)(c+H^2)\;\,\text{for any unit}\;\, X\in T_xM.
$$
Then the manifold $M^n$ is Einstein and $H$ is constant.
\vspace{1ex} 

If $H=0$ then $c>0$. Since $\Ric_M=(n-2)c$, it follows 
from the Theorem of Ejiri in \cite{E} that the submanifold is 
as in parts $(i)$ or $(ii)$ of the statement. \vspace{1ex}

We assume hereafter that $H>0$. We have to distinguish two cases
according to the dimension of the submanifold.
\vspace{1ex}

\noindent\emph{Case $n\geq 5$}.
Part $(ii)$ of Proposition \ref{prop} yields 
$k= n/2$ and that there are smooth Dupin
principal normal vector fields $\eta_1$ and $\eta_2$ of 
multiplicity $k$ and corresponding smooth distributions 
$E_1$ and $E_2$. Let $\{X_\ell\}_{1\leq\ell\leq n}$ be a 
smooth local orthonormal frame satisfying that
$E_{\eta_1}=\spa\left\{X_1,\dots,X_k\right\}$ and 
$E_{\eta_2}=\spa\left\{X_{k+1},\dots,X_n\right\}$. Then 
$\a_f(X_i,X_i)=\eta_1$ if $1\leq i\leq k$ and 
$\a_f(X_j,X_j)=\eta_2$ if $k+1\leq j\leq n$.
\vspace{1ex}

If follows from \eqref{ric} that 
$$
Ric_M(X)=(n-1)c\|X\|^2+n\<\mathcal{H},\a_f(X,X)\>-I\!I\!I(X)
\;\,\text{for any}\;\,X\in\mathcal X(M),
$$
where $I\!I\!I(X)=\sum_{\ell=1}^n\|\a_f(X,X_\ell)\|^2$ is 
the third fundamental form of $f$.
Since we have that $\mathcal H=(\eta_1+\eta_2)/2$, then
\be\label{H}
4H^2=\|\eta_1\|^2+\|\eta_2\|^2+2\<\eta_1,\eta_2\>.
\ee
Moreover, we have for $1\leq i\leq k$ that
$$
I\!I\!I(X_i)=\sum_{\ell=1}^n\|\a(X_\ell,X_i)\|^2
=\|\eta_1\|^2
$$
and
\begin{align*}
(n-2)(c+H^2)&=Ric_M(X_i)=(n-1)c+n\<\mathcal H,\a(X_i,X_i)\>
-I\!I\!I(X_i)\\
&=(n-1)c+k\<\eta_1+\eta_2,\eta_1\>-\|\eta_1\|^2.
\end{align*}
Thus
\be\label{ric1}
(n-2)(c+H^2)=(n-1)c+(k-1)\|\eta_1\|^2+k\<\eta_1,\eta_2\>.
\ee
Arguing similarly for $k+1\leq j\leq n$, we obtain
\be\label{ric2}
(n-2)(c+H^2)=(n-1)c+(k-1)\|\eta_2\|^2+k\<\eta_1,\eta_2\>.
\ee
It follows from \eqref{ric1} and \eqref{ric2} that
$\|\eta_1\|=\|\eta_2\|$, and hence \eqref{H} 
becomes
\be\label{H1}
2H^2=\|\eta_1\|^2+\<\eta_1,\eta_2\>.
\ee
Combining \eqref{ric1} with \eqref{H1} gives
\be\label{12}
c+\<\eta_1,\eta_2\>=0. 
\ee
Then, we conclude from \eqref{H1} that 
\be\label{1122}
\|\eta_1\|^2=\|\eta_2\|^2=2H^2+c.
\ee

The Codazzi equation for $f$ is easily seen to yield
\be\label{cod}
\<\nabla_XY,Z\>(\eta_i-\eta_j)=\<X,Y\>\nabla_Z^\perp\eta_i
\;\,\text{if}\;\,i\neq j
\ee
for all $X,Y\in E_i,Z\in E_j$.
Using \eqref{12} and \eqref{1122} then \eqref{cod} gives
$$
2\<\nabla_XY,Z\>H^2=\<X,Y\>\<\nabla_Z^\perp\eta_i,\eta_i\>=0
\;\,\text{for all}\;\,X,Y\in E_i\;\,\text{and}\;\,Z\in E_j,
\;i\neq j,
$$
that is, the distributions $E_1$ and $E_2$ are totally geodesic. 
Being simply connected, it is well-known that $M^n$ is a Riemannian 
product $M_1^k\times M_2^k$ (cf.\ Theorem $8.2$ in \cite{DT}) such that 
$TM_j^k=E_j$, $j=1,2$. Since the second fundamental form of $f$ 
is adapted to the distributions $E_1$ and $E_2$, then Theorem $8.4$ 
and Corollary $8.6$  in \cite{DT} imply that the submanifold is an 
extrinsic product of isometric immersions each of which is totally 
umbilical. Hence, if $c=0$ then the submanifold is a torus 
$\Sf^{n/2}(r/\sqrt{2})\times\Sf^{n/2}(r/\sqrt{2})$ in a sphere 
$\Sf^{n+1}(r)\subset \R^{n+2}$. If $c>0$, then the submanifold is a 
torus $\Sf^{n/2}(r/\sqrt{2})\times\Sf^{n/2}(r/\sqrt{2})$ in a sphere 
$\Sf^{n+1}(r)\subset\Sf^{n+2}(1/\sqrt{c})$.
\vspace{2ex}

\noindent\emph{Case $n=4$}. We have $k=2$ and $H_2(M^4;\Z)\neq 0$. 
Since $\Ric_M=2(c+H^2)$ then $\tau=8(c+H^2)$ and \eqref{tau} 
gives $S=4c+8H^2$. Thus equality at any point of the submanifold 
holds in the pinching condition $(1)$ in \cite{OV}. 
Since $\|\phi\|^2=S-4H^2$, it then follows from Proposition $16$ in 
\cite{OV} that the Bochner operator 
$\B^{[2]}\colon\Omega^2(M^4)\to\Omega^2(M^4)$, 
a certain symmetric endomorphism of the bundle of $2$-forms 
$\Omega^2(M^4)$, satisfies for any $\omega\in\Omega^2(M^4)$ 
the inequality
\be\label{ineq}
\<\B^{[2]}\omega,\omega\>
\geq (4c+8H^2-S)\|\omega\|^2=0.
\ee
 
We claim that there exists a nonzero $2$-form for which 
equality holds in \eqref{ineq} at any point. By the 
Bochner-Weitzenb\"ock formula the Laplacian of any 
$2$-form $\omega\in\Omega^2(M^4)$ is given by
$$
\Delta\omega=\n^*\n\omega+\B^{[2]}\omega,
$$
where $\n^*\n$ is the rough Laplacian. From this we obtain 
\be\label{boch.formula}
\<\Delta \omega,\omega\>
=\|\n\omega\|^2+\<\B^{[2]}\omega,\omega\>
+\frac{1}{2}\,\Delta \|\omega\|^2.
\ee 
If $\omega$ is an harmonic $2$-form, it follows from 
the maximum principle, \eqref{ineq} and \eqref{boch.formula}
and it is parallel. Then, for any harmonic $2$-form we have
that \eqref{ineq} holds as an equality at any point.
On the other hand, the Universal coefficient theorem of 
cohomology yields that the torsion subgroups of $H_1(M^4;\Z)$ 
and $H^2(M^4;\Z)$ are isomorphic 
(cf.\ \cite[p.\ 244 Corollary 4]{Sp}). 
Since $M^4$ is simply connected, we have that $H_1(M^n;\Z)=0$ 
and thus $H^2(M^4;\Z)$ is torsion free. Then the Poincar\'e 
duality yields that also $H_2(M^4;\Z)$ is torsion free. Hence, 
$0\neq H_2(M^4;\Z)=\Z^{\beta_2(M)}$. Thus $M^4$ supports a 
nonzero parallel harmonic $2$-form, and the claim has 
been proved.

From the claim and Proposition $16$ in \cite{OV} if follows
that the shape operator $A_\xi(x)$ at any $x\in M^4$ and 
for any $0\neq\xi\in N_fM(x)$
has at most two distinct eigenvalues with multiplicity $2$. 
We choose an orthonormal normal basis
$\{\xi_\alpha\}_{1\leq\alpha\leq m}$ at $x\in M^4$ such that the 
mean curvature vector is $\mathcal{H}(x)=H\xi_1$. By part $(ii)$ 
of Proposition \ref{prop} there exists an orthonormal basis 
$\{e_i\}_{1\leq i\leq 4}$  of $T_xM$ such that the corresponding 
shape operators $A_\a,1\leq\a\leq m$, are as
\be\label{A}
\begin{cases}
A_\a e_1=\rho_\a e_1+\kappa_\a e_3+\lambda_\a e_4\\
A_\a e_2=\rho_\a e_2+\mu_\a e_3+\nu_\a e_4\\
A_\a e_3=\kappa_\a e_1+\mu_\a e_2+\sigma_\a e_3\\
A_\a e_4=\lambda_\a e_1+\nu_\a e_2+\sigma_\a e_4,
\end{cases}
\ee
where
\be\label{1st}
\rho_1+\sigma_1=2H\;\;\text{and}\;\;\rho_\a+\sigma_\a=0
\;\;\text{for any}\;\;2\leq\a\leq m.
\ee
Now since each shape operator has at most two distinct 
eigenvalues with multiplicity $2$, it follows easily using 
\eqref{A} that 
\be\label{munu}
\nu_\a=\pm\kappa_\a\;\;\text{and}\;\;\mu_\a=\mp\lambda_\a
\;\;\text{for any}\;\;2\leq\a\leq m.
\ee

Since $M^4$ is Einstein, then
$\Ric_M(e_i,e_j)=0, i=1,2,\;j=3,4$. Using that
$$
\Ric_M(X,Y)=(n-1)c\<X,Y\>+\sum_{\a=1}^m(\mathrm{tr}A_\a)
\<A_\a X,Y\>-\sum_{\a=1}^m\<A^2_{\a}X,Y\>,
$$
together with \eqref{A} and \eqref{1st} the above yields 
that $\{e_i\}_{1\leq i\leq 4}$ diagonalizes $A_1$. 
Moreover, from $\Ric_M(e_j)=2(c+H^2), 1\leq j\leq 4$,
\eqref{ric} and \eqref{A} we obtain that
$$
2H^2-c=
\begin{cases}
4H\rho_1-\sum_{\a\geq 1}\rho^2_\a-\|\a_{13}\|^2-\|\a_{14}\|^2\\
4H\rho_1-\sum_{\a\geq 1}\rho^2_\a-\|\a_{23}\|^2-\|\a_{24}\|^2\\
4H\sigma_1-\sum_{\a\geq 1}\rho^2_\a-\|\a_{13}\|^2-\|\a_{23}\|^2\\
4H\sigma_1-\sum_{\a\geq 1}\rho^2_\a-\|\a_{14}\|^2-\|\a_{24}\|^2.
\end{cases}
$$
This implies that $\|\a_{13}\|
=\|\a_{24}\|,\|\a_{23}\|=\|\a_{14}\|$ 
and then that $\rho_1=\sigma_1=H$, namely, the submanifold is 
pseudo-umbilical. 

Being the submanifold pseudo-umbilical, we have from \eqref{A} 
and \eqref{munu} that
\be\label{alla} 
\a_{11}=\a_{22},\;\a_{33}=\a_{44},\; \a_{12}=\a_{34}=0,
\;\a_{13}=\pm\a_{24},\; \a_{23}=\mp\a_{14}.
\ee
Thus the vector subspace $N_1(x)\subset N_fM(x)$ spanned by 
the second fundamental form at $x\in M^n$ satisfies
$\dim N_1(x)\leq 4$ at any $x\in M^n$. 

We claim that the mean curvature vector field is parallel 
in the normal bundle.
Let $U$ be an open subset of $M^n$ were the subspaces $N_1(x)$ 
have constant dimension $r$, $1\leq r\leq 4$, and hence $N_1|_U$
is a smooth vector subbundle of the normal bundle. Then let 
$\{e_i\}_{1\leq i\leq 4}$ be a local smooth orthonormal frame 
with respect to which the second fundamental form is as in 
\eqref{alla}. From the Codazzi equation
$$
(\n_{e_i}^\perp\a)(e_j, e_k)=(\n_{e_j}^\perp\a)(e_i, e_k),
$$
we obtain that $\n_{e_i}^\perp\a(e_j,e_k)\in N_1|_U$ for any 
$1\leq i,j,k\leq 4$. Hence $N_1|_U$ is a parallel subbundle 
of the normal bundle and, consequently, $f|_U$ reduces its 
codimension, that is, it is a composition $f|_U=i\circ g$ 
where $i\colon\Q^{4+r}_c\to\Q^{4+m}_c$ is a totally geodesic 
inclusion and $g$ is an isometric immersion into $\Q^{4+r}_c$.

Since the submanifold is pseudo-umbilical, 
from the Codazzi equation
$$
(\n_XA_1)Y-(\n_YA_1)X=A_{\n^\perp_X\xi_1}Y-A_{\n^\perp_Y\xi_1}X,
\;\,\text{for all}\;\,X,Y\in\mathcal X(M),
$$
we obtain that $\n^\perp_X\xi_1\in (N_1|_U)^\perp$ for any
$X\in\mathcal X(M)$.
Hence the mean curvature vector field is parallel in the 
normal bundle of $f$ along any open subset where the dimension 
of the first normal space is constant. By continuity this
is the case globally. Thus, it is an elementary fact that 
the submanifold decomposes as $f=j\circ g$, where 
$g\colon M^n\to\Q^{4+p}_{\tilde c}$ is a minimal submanifold 
and $j\colon\Q^{4+p}_{\tilde c}\to\Q^{4+m}_c$ is totally 
umbilical with $\tilde c=c+H^2$. It now follows from the 
result of Ejiri \cite{E} that the submanifold is as in 
parts $(i)$ or $(ii)$ of the statement of the theorem.
\vspace{2ex}\qed

\noindent Marcos Dajczer is  partially supported 
by the grant PID2021-124157NB-I00 funded by 
MCIN/AEI/10.13039/501100011033/ `ERDF A way of making Europe',
Spain, and are also supported by Comunidad Aut\'{o}noma de la 
Regi\'{o}n de Murcia, Spain, within the framework of the Regional 
Programme in Promotion of the Scientific and Technical Research 
(Action Plan 2022), by Fundaci\'{o}n S\'{e}neca, Regional Agency 
of Science and Technology, REF, 21899/PI/22.

Theodoros Vlachos thanks the Department of Mathematics 
of the University of Murcia where part of this work was 
done for its cordial hospitality during his visit. 
He was supported by the grant PID2021-124157NB-I00 
funded by MCIN/AEI/10.13039/501100011033/ `ERDF 
A way of making Europe', Spain.

\noindent Marcos Dajczer\\
Departamento de Matemáticas\\ 
Universidad de Murcia, Campus de Espinardo\\ 
E-30100 Espinardo, Murcia, Spain\\
e-mail: marcos@impa.br
\bigskip

\noindent Theodoros Vlachos\\
University of Ioannina \\
Department of Mathematics\\
Ioannina -- Greece\\
e-mail: tvlachos@uoi.gr
\end{document}